\documentclass[10pt]{amsart}

\newcounter{minutes}
\divide\time by 60
\newcounter{hours}
\multiply\time by 60 \addtocounter{minutes}{-\time}

\usepackage{amsfonts}
\usepackage{amsthm}
\usepackage{amsmath}
\usepackage{graphicx}
\usepackage{amscd,amssymb,amsthm}
\usepackage[mathscr]{eucal}

\newtheorem*{lemma}{Lemma}
\newtheorem{theorem}{Theorem}
\newtheorem*{remark}{Remark}
\newtheorem{corollary}{Corollary}[theorem]
\newtheorem*{conjecture}{Conjecture}

\title{Tur\'an determinants of Bessel functions}

\author[]{\'Arp\'ad Baricz}
\address{Department of Economics, Babe\c{s}-Bolyai University,
Cluj-Napoca 400591, Romania} \email{bariczocsi@yahoo.com}

\author[\'A. Baricz, T.K. Pog\'any]{Tibor K. Pog\'any}
\address{Faculty of Maritime Studies, University of Rijeka, Rijeka 51000, Croatia}
\email{poganj@brod.pfri.hr}



\begin{document}

\begin{center}
\texttt{File:~\jobname .tex,
          printed: \number\year-0\number\month-\number\day,
          \thehours.\ifnum\theminutes<10{0}\fi\theminutes}
\end{center}

\maketitle

\begin{center}
{Dedicated to Agata, Bor\'oka and Kopp\'any}
\end{center}

\begin{abstract}
In this paper first we survey the Tur\'an type inequalities and
related problems for the Bessel functions of the first kind. Then we
extend the known higher order Tur\'an type inequalities for Bessel
functions of the first kind to real parameters and we deduce new
closed integral representation formulae for the second kind Neumann
type series of Bessel functions of the first kind occurring in the
study of Tur\'an determinants of Bessel functions of the first kind.
At the end of the paper we prove a Tur\'an type inequality for the
Bessel functions of the second kind.
\end{abstract}

\section{\bf Introduction}

Consider the Bessel function of the first kind and order $\nu$,
defined by \cite[p. 40]{watson}
$$
J_{\nu}(x)=\sum_{n\geq0}\frac{(-1)^n}{n!\Gamma(n+\nu+1)}\left(\frac{x}{2}\right)^{2n+\nu},
$$
which is a particular solution of the second-order homogeneous
Bessel differential equation \cite[p. 38]{watson}
\begin{equation}\label{eqbessel}
x^2y''(x)+xy'(x)+(x^2-\nu^2)y(x)=0.
\end{equation}
In this paper our aim is to survey the Tur\'an type and related
inequalities for the Bessel functions of the first kind. The Tur\'an
type inequalities have a long history and because of some
applications, recently many authors have deduced new Tur\'an type
inequalities for some special functions (e.g. Kr\"atzel functions,
Gauss and Kummer hypergeometric functions, generalized
hypergeometric functions, modified Bessel functions, Euler gamma
function). For more details we refer to the papers [3--12], \cite{bepe, barnard, dimitrov, dimitrov2, karp, laforgiaism, segura} and
to the references contained therein. Motivated by the vivid interest
on Tur\'an type inequalities in this paper we would like to show how
a simple and beautiful inequality, like Tur\'an inequality, can
occur in many problems of analysis, especially classical analysis.
This survey treats interesting and nice results of different
authors.

The paper is organized as follows: in section 2 we survey the
Tur\'an type inequalities for Bessel functions of the first kind.
Here Tur\'an and Hankel determinants of Bessel functions of the
first kind, generalized Tur\'an expression of Bessel functions,
properties of the relative minima and maxima of the Tur\'an
expression, log-concavity of Bessel functions of the first kind with
respect to their order, zeros of the Tur\'an expression in the
complex plane and lower bounds for the logarithmic derivative of
Bessel functions of the first kind are discussed. Section 2 is
closed with a conjecture for Bessel functions of the first kind. In
section 3 we will show some higher order Tur\'an type inequalities
for Bessel functions of the first kind. Relevant connections of the
results of section 3 and section 2 are also pointed out. In section
4 we deduce a closed integral representation formula for second kind
Neumann type series of Bessel functions, which occur in section 2.
Finally, in section 5 we present a Tur\'an type inequality for the
Bessel function of the second kind, solving a recent open problem
from \cite{bams}, and we apply this result to obtain a lower bound
for the logarithmic derivative of Bessel functions of the second
kind.

\section{\bf Tur\'an type inequalities for Bessel functions of the first kind and related problems}
\setcounter{equation}{0}

In this section we are mainly interested on Tur\'an type
inequalities for Bessel functions of the first kind and some related
problems. The results are presented in the chronological order and
relevant connections between them are pointed out.

\subsection{Tur\'an type inequalities for Bessel functions of the first kind}
Let us start with von Lommel's formula \cite[p. 152]{watson} for
Bessel functions of the first kind
   \begin{equation} \label{lommel}
      x^2\left[J_{\nu}^2(x)-J_{\nu-1}(x)J_{\nu+1}(x)\right] = \sum_{n\geq0}(\nu+1+2n)J_{\nu+1+2n}^2(x),
   \end{equation}
which is valid for all $x\in\mathbb{R}$ and $\nu>-1.$ Then clearly
we obtain the following Tur\'an type inequality for all
$x\in\mathbb{R}$ and $\nu>-1$
\begin{equation}\label{turan1}
\Delta_{\nu}(x)=J_{\nu}^2(x)-J_{\nu-1}(x)J_{\nu+1}(x)\geq0.
\end{equation}
The above simple and beautiful inequality had attracted many
mathematicians and has been proved in several ways. In 1950 and 1951
Sz\'asz \cite{szasz1, szasz2} deduced \eqref{turan1} and its sharper
form
\begin{equation}\label{turan2}
J_{\nu}^2(x)-J_{\nu-1}(x)J_{\nu+1}(x)\geq\frac{1}{\nu+1}J_{\nu}^2(x),
\end{equation}
where $x\in\mathbb{R}$ and $\nu>0,$ by using the recurrence
relations for the Bessel functions of the first kind. We note that,
since
\begin{equation}\label{limit}\lim_{x\to0}\frac{\Delta_{\nu}(x)}{J_{\nu}^2(x)}=
1-\lim_{x\to0}\frac{J_{\nu-1}(x)J_{\nu+1}(x)}{J_{\nu}^2(x)}=\frac{1}{\nu+1},\end{equation}
the constant $1/(\nu+1)$ on the right-hand side of \eqref{turan2} is
the best possible. See \cite{bams} or the last subsection of this
section for more details. In 1951 Thiruvenkatachar and Nanjundiah
\cite{thiru} proved also \eqref{turan1} by using recurrence
relations and deduced \eqref{turan2} by means of the following
formula
\begin{equation}\label{turan3}
\nu
J_{\nu}^2(x)-(\nu+1)J_{\nu-1}(x)J_{\nu+1}(x)=J_{\nu+1}^2(x)+\nu\left[J_{\nu+1}^2(x)-J_{\nu}(x)J_{\nu+2}(x)\right].
\end{equation}
It is interesting to note here that the sharp inequality
\eqref{turan2} in view of \eqref{turan3} is actually a consequence
of the weaker result \eqref{turan1}. With other words, for $\nu>0$
the Tur\'an type inequalities \eqref{turan1} and \eqref{turan2} are
equivalent. Note also that by using \eqref{turan3} Thiruvenkatachar
and Nanjundiah \cite{thiru} proved the following formula
   \begin{equation} \label{turan4}
      \Delta_{\nu}(x)=\frac{1}{\nu+1}J_{\nu}^2(x)+\frac{2}{\nu+2}J_{\nu+1}^2(x)
      +2\nu\sum_{n\geq2}\frac{J_{\nu+n}^2(x)}{(\nu+n-1)(\nu+n+1)}
   \end{equation}
which clearly implies the sharp Tur\'an type inequality
\eqref{turan2}. Moreover, it should be mentioned that in view of
\eqref{turan4} the Tur\'an type inequality \eqref{turan2} can be
improved as follows
\begin{equation}\label{turan5}J_{\nu}^2(x)-J_{\nu-1}(x)J_{\nu+1}(x)\geq
\frac{1}{\nu+1}J_{\nu}^2(x)+\frac{2}{\nu+2}J_{\nu+1}^2(x),\end{equation}
where $x\in\mathbb{R}$ and $\nu\geq0.$ Observe that the Tur\'an type
inequality \eqref{turan5} is sharp as $x=0$ or $\nu=0.$ This can be
seen easily by using that
$$\lim_{x\to0}\left[2\nu\sum_{n\geq2}\frac{J_{\nu+n}^2(x)}{(\nu+n-1)(\nu+n+1)}\right]=0$$
and by noting that $J_{-1}(x)=J_1(x)$ for all real $x.$

Another proof for the Tur\'an type inequality \eqref{turan1} was
given by Skovgaard \cite{skov} in 1954. Skovgaard's proof is
somewhat similar to the proof of Thiruvenkatachar and Nanjundiah
\cite{thiru}, however in \cite{skov} there is used in addition the
infinite product representation \cite[p. 498]{watson} of the Bessel
function of the first kind
\begin{equation}\label{weier}2^{\nu}\Gamma(\nu+1)x^{-\nu}J_{\nu}(x)=\prod_{n\geq1}\left(1-\frac{x^2}{j_{\nu,n}^2}\right).\end{equation}

\subsection{Tur\'an and Hankel determinants of Bessel functions of the first kind}
A generalization of the Tur\'an type inequality \eqref{turan1} has
been proved by Karlin and Szeg\H o in their mammoth work \cite[p.
130]{karlin} in 1960. They proved that for all $x>0,$ $\nu>-1$ and
$n$ even natural number the following inequalities are valid for the
Tur\'an and Hankel determinants of Bessel functions
$$(-1)^{\frac{n}{2}}\left|\begin{array}{cccc}\widetilde{J}_{\nu}(x) & \widetilde{J}_{\nu+1}(x) & \dots & \widetilde{J}_{\nu+n-1}(x)\\
\widetilde{J}_{\nu+1}(x) & \widetilde{J}_{\nu+2}(x) & \dots &
\widetilde{J}_{\nu+n}(x)\\ \vdots & \vdots &  &
\vdots\\\widetilde{J}_{\nu+n-1}(x) & \widetilde{J}_{\nu+n}(x) &
\dots & \widetilde{J}_{\nu+2n-2}(x)\end{array}\right|\geq0$$ and
$$(-1)^{\frac{n}{2}}\left|\begin{array}{cccc}\widetilde{J}_{\nu}(x) & \widetilde{J}_{\nu}'(x) & \dots & \widetilde{J}_{\nu}^{(n-1)}(x)\\
\widetilde{J}_{\nu}'(x) & \widetilde{J}_{\nu}''(x) & \dots &
\widetilde{J}_{\nu}^{(n)}(x)\\ \vdots & \vdots &  &
\vdots\\\widetilde{J}_{\nu}^{(n-1)}(x) &
\widetilde{J}_{\nu}^{(n)}(x) & \dots &
\widetilde{J}_{\nu}^{(2n-2)}(x)\end{array}\right|\geq0,$$ where
$\widetilde{J}_{\nu}(x)=x^{-\frac{\nu}{2}}J_{\nu}(2\sqrt{x}).$
Moreover, Karlin and Szeg\H o \cite[p. 131]{karlin} pointed out that
in the above Tur\'an determinant we can cancel out the factors of
powers of $x$ to have
$$(-1)^{\frac{n}{2}}\left|\begin{array}{cccc}{J}_{\nu}(2\sqrt{x}) & {J}_{\nu+1}(2\sqrt x) & \dots & {J}_{\nu+n-1}(2\sqrt x)\\
{J}_{\nu+1}(2\sqrt x) & {J}_{\nu+2}(2\sqrt x) & \dots &
{J}_{\nu+n}(2\sqrt x)\\ \vdots & \vdots & &
\vdots\\{J}_{\nu+n-1}(2\sqrt x) & {J}_{\nu+n}(2\sqrt x) & \dots &
{J}_{\nu+2n-2}(2\sqrt x)\end{array}\right|\geq0,$$ which can be
rewritten also as
$$(-1)^{\frac{n}{2}}\left|\begin{array}{cccc}{J}_{\nu}(x) & {J}_{\nu+1}(x) & \dots & {J}_{\nu+n-1}(x)\\
{J}_{\nu+1}(x) & {J}_{\nu+2}(x) & \dots & {J}_{\nu+n}(x)\\
\vdots & \vdots & & \vdots\\{J}_{\nu+n-1}(x) & {J}_{\nu+n}(x) &
\dots & {J}_{\nu+2n-2}(x)\end{array}\right|\geq0.$$

\subsection{Generalized Tur\'an expression of Bessel functions}
Now, we are going to present a genera\-lization of von Lommel's result
\eqref{lommel}. In 1961 Al-Salam \cite{salam} proved the following
result:
   \begin{align} \label{salam} \nonumber
      \frac{4^m(2m)!}{x^{2m}m!(m-1)!} &\sum_{k\geq0}(\nu+m+2k){(k+1)}_{m-1}{(\nu+k+1)}_{m-1}J_{\nu+m+2k}^2(x)\\
                                      &=\sum_{n=-m}^{m}(-1)^n\left(\begin{array}{c}2m\\
                                        m-n\end{array}\right)J_{\nu-n}(x)J_{\nu+n}(x),
   \end{align}
where ${(a)}_0=1$ for $a\neq0$ and
${(a)}_n=a(a+1)(a+2)\dots(a+n-1)=\Gamma(a+n)/\Gamma(a)$ for
$n\in\{1,2,\dots\}$ is the well-known Pochhammer (or Appell) symbol
in terms of the Euler gamma function. For reader's convenience we
note that there is a small typographical mistake in the original
formula of Al-Salam: the term ${(n+k-1)}_{m-1}$ in \cite[Eq.
(1.2)]{salam} should be written as ${(n+k+1)}_{m-1}.$ We also note
that in \cite{salam} the author does not discussed the range of
validity of \eqref{salam}, however, a close inspection of the proof
of \eqref{salam} reveals that this formula holds for all
$m\in\{0,1,\dots\},$ $x\in\mathbb{R}$ and $\nu>-1.$ Furthermore,
\eqref{salam} implies the extended Tur\'an type inequality
\begin{equation}\label{salamt}
\sum_{n=-m}^{m}(-1)^n\left(\begin{array}{c}2m\\
m-n\end{array}\right)J_{\nu-n}(x)J_{\nu+n}(x)\geq0,
\end{equation}
where $m\in\{0,1,\dots\},$ $x\in\mathbb{R}$ and $\nu>-1.$ Observe
that when $m=1$ in \eqref{salam} and \eqref{salamt}, then we
reobtain \eqref{lommel} and \eqref{turan1}. Moreover, if we choose
$m=2$ in \eqref{salam}, then we obtain
\begin{align*}192x^{-4}\sum_{k\geq0}&(\nu+2k+2)(k+1)(\nu+k+1)J_{\nu+2k+2}^2(x)\\
&=6J_{\nu}^2(x)-8J_{\nu-1}(x)J_{\nu+1}(x)+2J_{\nu-2}(x)J_{\nu+2}(x)\\
&=6\Delta_{\nu}(x)+2\left[J_{\nu-2}(x)J_{\nu+2}(x)-J_{\nu-1}(x)J_{\nu+1}(x)\right],
\end{align*}
which immediately implies the new Tur\'an type inequality
$$3\Delta_{\nu}(x)\geq J_{\nu-1}(x)J_{\nu+1}(x)-J_{\nu-2}(x)J_{\nu+2}(x),$$
where $x\in\mathbb{R}$ and $\nu>-1.$

\subsection{Properties of the relative minima and maxima of Tur\'an
expression} In 1961 Lakshmana Rao \cite{rao} proved the following
result concerning the Tur\'an expression $\Delta_{\nu}(x)$:
\begin{theorem}\label{thlak}
Let $\nu>0$ and $x\in\mathbb{R}.$ The relative maxima (denoted by
$M_{\nu,n}$) of the function $x\mapsto\Delta_{\nu}(x)$ occur at the
zeros of $J_{\nu-1}$ and the relative minima (denoted by
$m_{\nu,n}$) occur at the zeros of $J_{\nu+1}.$ The values
$M_{\nu,n}$ and $m_{\nu,n}$ can be expressed as
$$M_{\nu,n}=\Delta_{\nu}(j_{\nu-1,n})=J_{\nu}^2(j_{\nu-1,n})$$
and
$$m_{\nu,n}=\Delta_{\nu}(j_{\nu+1,n})=J_{\nu}^2(j_{\nu+1,n}).$$
Since the quantities $M_{\nu,n}$ and $m_{\nu,n}$ are positive for
all values of $\nu,$ we have the Tur\'an type inequality
\eqref{turan1}. Moreover, the sequences $\{M_{\nu,n}\}_{n\geq 1}$
and $\{m_{\nu,n}\}_{n\geq1}$ are decreasing for $\nu$ sufficiently
large and in addition, for each $\nu$ we have $M_{\nu,n}>m_{\nu,n}.$
Finally, for a fixed value of $n,$ the functions $\nu\mapsto
M_{\nu,n}$ and $\nu\mapsto m_{\nu,n}$ are decreasing.
\end{theorem}

\subsection{Log-concavity of Bessel functions of the first kind with respect to their order}
In 1978 Ismail and Muldoon \cite[Lemma 2.3]{im}, by using von
Neumann's formula for product of two Bessel functions with different
order, proved that for each fixed $\beta,$ where $0<\beta\leq1,$ and
$x\neq j_{\nu,n},$ where $n\in\{1,2,\dots\},$ the function $\nu
\mapsto J_{\nu+\beta}(x)/J_{\nu}(x)$ is decreasing on
$[-(\beta+1)/2,\infty).$ As it was pointed out later in 1997 by
Muldoon \cite{mul}, this actually implies that for fixed $x>0$ and
$x\neq j_{\nu,n},$ $n\in\{1,2,\dots\},$ the function $\nu\mapsto
J_{\nu}(x)$ is logarithmically concave on $(-1,\infty).$ And this
clearly yields the Tur\'an type inequality \eqref{turan1} for
$\nu>0$ and $x>0$ such that $x\neq j_{\nu,n},$ $n\in\{1,2,\dots\}.$
We note that the remained part of the Tur\'an type inequality
\eqref{turan1} when $x=j_{\nu,n},$ $n\in\{1,2,\dots\},$ can be
verified easily by using the recurrence relations
\begin{equation}\label{rec1}xJ_{\nu}'(x)+\nu J_{\nu}(x)=xJ_{\nu-1}(x)\end{equation} and
\begin{equation}\label{rec2}xJ_{\nu}'(x)-\nu J_{\nu}(x)=-xJ_{\nu+1}(x).\end{equation}
More precisely, in view of the formulas \eqref{rec1} and
\eqref{rec2} for all $n\in\{1,2,\dots\}$ and all $\nu>-1$ we have
$$\Delta_{\nu}(j_{\nu,n})=-J_{\nu-1}(j_{\nu,n})J_{\nu+1}(j_{\nu,n})=\left[J_{\nu}'(j_{\nu,n})\right]^2>0.$$

We note that some of the results (and proofs) of Sz\'asz
\cite{szasz2} concerning the Tur\'an type inequalities
\eqref{turan1} and \eqref{turan2} have been reobtained later in 1991
by Joshi and Bissu \cite{JB}. Moreover, it should be mentioned that
in 1997 Bustoz and Ismail \cite{bustoz} proved also the Tur\'an type
inequality \eqref{turan2} by showing the positivity of Tur\'an
determinants of symmetric Pollaczek polynomials, von Lommel
polynomials and $q-$Bessel functions. Finally, for applications of
the Tur\'an type inequality \eqref{turan2} and \cite[Eq. 2.9]{JB} in
the extension of some trigonometric inequalities (like
Mitrinovi\'c-Adamovi\'c and Wilker) to Bessel functions we refer to
the papers \cite{bexpo} and \cite{bsandor}.

\subsection{Zeros of the Tur\'an expression $\Delta_{\nu}$ in the complex plane}
In 2001 Kravanja and Verlinden \cite{kra} applied the Tur\'an type
inequality \eqref{turan1} in the study of the zeros of
$\Delta_{\nu}.$ More precisely, the Tur\'an-type inequality
\eqref{turan1} implies that no zero exists on the positive real
axis, and, by symmetry, neither on the negative real axis. Moreover,
Kravanja and Verlinden \cite{kra} proved the following result:
\begin{theorem}
Except for $z=0,$ the function $z\mapsto \Delta_{\nu}(z)$ has no
zeros on the imaginary axis. Moreover, all the zeros of $z\mapsto
\Delta_{\nu}(z)$ that lies in $\mathbb{C}\setminus\{0\}$ are simple.
\end{theorem}

The zeros of the function $z\mapsto \Delta_{\nu}(z)$ play an
important role in certain physical applications. For example,
MacDonald \cite{mac} used the zeros of $\Delta_1$ to plot
representative streamlines for the steady motion of a viscous fluid
in a long tube, of constant radius, which rotates about its axis
with an angular velocity that changes discontinuously at $z=0$ from
one constant value to another of the same sign.

It is also worth mentioning here that the Tur\'an expression
$\Delta_{\nu}(x)$ appears in other physical applications. More
precisely, the series
$$\sum_{n=-\infty}^{\infty}\left[\Delta_n(x)\right]^2$$
and
$$\sum_{m=-\infty}^{\infty}\sum_{n=-\infty}^{\infty}\Delta_n(x)\Delta_m(x)\Delta_{n-m}(\mu
x),$$ where $\mu$ is a constant, appear in the investigation into
multiple scattering of acoustic waves by random configurations of
penetrable circular cylinders. For more details see \cite{martin},
where the asymptotic behavior for large positive $x$ of the above
series were investigated.

\subsection{Lower bounds for the logarithmic derivative of Bessel
functions of the first kind} In this subsection our aim is to deduce
some new inequalities for the logarithmic derivative of $J_{\nu}$
and for the quotient $J_{\nu}/J_{\nu-1}$ by using the Tur\'an type
inequalities \eqref{turan1} and \eqref{turan2}. Observe that
combining \eqref{turan1} with \eqref{rec1} and \eqref{rec2} we
immediately get the inequality
\begin{equation}\label{bound1}
\left[xJ_{\nu}'(x)\right]^2+(x^2-\nu^2)J_{\nu}^2(x)\geq0,
\end{equation}
which holds for all $\nu>-1$ and $x\in\mathbb{R}.$ Moreover,
combining \eqref{bound1} with \eqref{rec1} we arrive at
\begin{equation}\label{bound2}
xJ_{\nu}^2(x)-2\nu J_{\nu-1}(x)J_{\nu}(x)+xJ_{\nu-1}^2(x)\geq0,
\end{equation}
where $\nu>-1$ and $x\in\mathbb{R}.$ The inequality \eqref{bound1}
implies the inequality
\begin{equation}\label{bound3}
\left[\frac{xJ_{\nu}'(x)}{J_{\nu}(x)}\right]^2\geq \nu^2-x^2,
\end{equation}
where $\nu>-1$ and $x\neq j_{\nu,n},$ $n\in\{1,2,\dots\}.$ While
\eqref{bound2} implies the following
\begin{equation}\label{bound4}
\frac{J_{\nu}(x)}{J_{\nu-1}(x)}\geq \frac{\nu+\sqrt{\nu^2-x^2}}{x},
\end{equation}
where $\nu>0,$ $x\in(0,\nu]$ such that $x\neq j_{\nu-1,n},$
$n\in\{1,2,\dots\}.$

Now, we are going to improve the above inequalities. Let
$\mu=\nu/(\nu+1).$ Then combining \eqref{turan2} with \eqref{rec1}
and \eqref{rec2} we immediately get the inequality
\begin{equation}\label{boundn1}
\left[xJ_{\nu}'(x)\right]^2+\left(\mu
x^2-\nu^2\right)J_{\nu}^2(x)\geq0,
\end{equation}
which holds for all $\nu>-1$ and $x\in\mathbb{R}.$ Moreover,
combining \eqref{boundn1} with \eqref{rec1} we obtain the following
\begin{equation}\label{boundn2}
\mu xJ_{\nu}^2(x)-2\nu J_{\nu-1}(x)J_{\nu}(x)+xJ_{\nu-1}^2(x)\geq0,
\end{equation}
where $\nu>-1$ and $x\in\mathbb{R}.$ The inequality \eqref{boundn1}
implies the inequality
\begin{equation}\label{boundn3}
\left[\frac{xJ_{\nu}'(x)}{J_{\nu}(x)}\right]^2\geq \nu^2-\mu
x^2,\end{equation} where $\nu>-1$ and $x\neq j_{\nu,n},$
$n\in\{1,2,\dots\}.$ While \eqref{boundn2} implies the following
\begin{equation}\label{boundn4}
\frac{J_{\nu}(x)}{J_{\nu-1}(x)}\geq \frac{\nu+\sqrt{\nu^2-\mu
x^2}}{\mu x},\end{equation} where $\nu>0,$
$x\in(0,\sqrt{\nu(\nu+1)}]$ such that $x\neq j_{\nu-1,n},$
$n\in\{1,2,\dots\}.$

Observe that since $\nu>-1$ we have $\mu<1$ and consequently
$$\nu^2-\mu
x^2>\nu^2-x^2$$ and
$$\frac{\nu+\sqrt{\nu^2-\mu
x^2}}{\mu x}>\frac{\nu+\sqrt{\nu^2-\mu
x^2}}{x}>\frac{\nu+\sqrt{\nu^2-x^2}}{x}$$ hold for all $\nu,x>0.$
These in turn imply that \eqref{boundn3} improves \eqref{bound3} and
\eqref{boundn4} improves \eqref{bound4}.

Finally, it is important to note here that similar inequalities to
those given in \eqref{bound1} -- \eqref{boundn4} for modified Bessel
functions of the first and second kinds $I_{\nu}$ and $K_{\nu}$ were
given recently by Baricz \cite{bpams}, Laforgia and Natalini
\cite{lana}, and Segura \cite{segura}. For more details see also
\cite{bepe} and the references therein.

\subsection{On a conjecture for Bessel functions of the first kind}
In this subsection we recall a recent conjecture of the first author
\cite{bams}. Consider the function
$\Phi:(0,\infty)\setminus\Xi\to(0,\infty),$ defined by (see
\cite{bams})
$$\Phi_{\nu}(x)=\frac{\Delta_{\nu}(x)}{J_{\nu}^2(x)}=1-\frac{J_{\nu-1}(x)J_{\nu+1}(x)}{J_{\nu}^2(x)}=\sum_{n\geq
1}\frac{4j_{\nu,n}^2}{(x^2-j_{\nu,n}^2)^2},$$ where
$\Xi=\{j_{\nu,n}:n\geq 1\}$ is the set of the zeros of the Bessel
function of the first kind of order $\nu.$ Based on numerical
experiments we conjecture (see \cite{bams}) that Sz\'asz's result
\eqref{turan2} can be improved as follows.

\begin{conjecture}
Let $\nu>0.$ Then the equation $\Phi'_{\nu}(x)=0$ has infinitely
many roots. Denoting with $\alpha_{\nu,n}$ these roots, where
$\alpha_{\nu,n}\in (j_{\nu,n},j_{\nu,n+1})$ for all
$n\in\{1,2,\dots\},$ the Tur\'an type inequality \eqref{turan2} can
be improved as follows
$$J_{\nu}^2(x)-J_{\nu-1}(x)J_{\nu+1}(x)>\beta_{\nu,n}J_{\nu}^2(x),$$
where $x\in(j_{\nu,n},j_{\nu,n+1}),$
$\beta_{\nu,n}=\Phi_{\nu}(\alpha_{\nu,n})$ for all
$n\in\{0,1,\dots\}$ and the sequence $\{\beta_{\nu,n}\}_{n\geq0}$ is
strictly increasing, where $\alpha_{\nu,0}=j_{\nu,0}=0$ and
$\beta_{\nu,0}=\Phi(\alpha_{\nu,0}^+)=1/(\nu+1).$
\end{conjecture}

Let $n\in\{1,2,{\dots}\}$ be fixed. Since
$$\Phi'(j_{\nu-1,n})=\frac{2\nu}{j_{\nu-1,n}}>0\ \ \mbox{and}\ \ \Phi'(j_{\nu+1,n})=-\frac{2\nu}{j_{\nu+1,n}}<0,$$ we obtain that
$\alpha_{\nu,n}\in(j_{\nu+1,n},j_{\nu-1,n+1}).$ Now, since
(according to Theorem \ref{thlak}) the function $x\mapsto
\Delta_{\nu}(x)$ is increasing on $(j_{\nu+1,n},j_{\nu-1,n+1}),$ we
obtain that
\begin{align*}\Delta_{\nu}'(\alpha_{\nu,n})&=({2}/{\alpha_{\nu,n}})J_{\nu-1}(\alpha_{\nu,n})J_{\nu+1}(\alpha_{\nu,n})\\&=
({2}/{\alpha_{\nu,n}})J_{\nu}^2(\alpha_{\nu,n})\left(1-\Phi_{\nu}(\alpha_{\nu,n})\right)>0,\end{align*}
which in turn implies that
$\beta_{\nu,n}=\Phi_{\nu}(\alpha_{\nu,n})<1$ for all $\nu>0$ and
$n\in\{1,2,{\dots}\}.$ It is also worth to mention here that the
relative minima $\beta_{\nu,n}$ does not occur in the zeros
$j_{\nu,n+1}'$ of $J_{\nu}',$ although it seems from graphics. For
example, for $\nu=2$ we have $\alpha_{2,1}=6.690090363,$
$\alpha_{2,2}=9.965082278,$ $\alpha_{2,3}=13.1685359,$
$\alpha_{2,4}=16.3465786$ and $j_{2,2}'=6.706133194,$
$j_{2,3}'=9.969467823,$ $j_{2,4}'=13.1703708,$
$j_{2,5}'=16.3475223.$

\section{\bf Higher order Tur\'an type inequalities for Bessel functions of the first kind}
\setcounter{equation}{0}

This section is devoted to the study of higher order Tur\'an type
inequalities and related problems for the Bessel functions of the
first kind. The results stated here complement naturally those
presented in the previous section.

\subsection{The Laguerre-P\'olya class of entire functions and
necessary conditions} By definition an entire function is a function
of a complex variable which is holomorphic in the entire complex
plane and can be represented by an everywhere convergent power
series
$$f(z)=a_0+a_1z+a_2z^2+a_3z^3+{\dots}+a_nz^n+{\dots}.$$
An entire function $f$ is called real if it is real on the real
axis, or equivalently, if it has only real coefficients in its power
series. Now, by definition the real entire function $\varphi,$
defined by
\begin{equation}\label{entire}\varphi(x)=\varphi(x;t)=\sum_{n\geq0}b_n(t)\frac{x^n}{n!},\end{equation} is said
to be in the Laguerre-P\'olya class (denoted by $\mathcal{LP}$), if
$\varphi(x)$ can be expressed in the form
$$\varphi(x)=cx^de^{-\alpha x^2+\beta x}\prod_{i=1}^{\omega}\left(1-\frac{x}{x_i}\right)e^{\frac{x}{x_i}},\ \ \ 0\leq\omega\leq\infty,$$
where $c$ and $\beta$ are real, $x_i$'s are real and nonzero for all
$i\in\{1,2,{\dots},\omega\},$ $\alpha\geq0,$ $d$ is a nonnegative
integer and $\sum_{i=1}^{\omega}x_i^{-2}<\infty.$ If $\omega=0,$
then, by convention, the product is defined to be $1.$ For the
various properties of the functions in the Laguerre-P\'olya class we
refer to \cite{craven1,craven2,norfolk,csordas} and to the
references contained therein. We note that in fact a real entire
function $\varphi$ is in the Laguerre-P\'olya class if and only if
$\varphi$ can be uniformly approximated on disks around the origin
by a sequence of polynomials with only real zeros. This in turn
implies that the class $\mathcal{LP}$ is closed under
differentiation, that is, if $\varphi\in\mathcal{LP}$, then
$\varphi^{(n)}\in\mathcal{LP}$ for all $n$ nonnegative integer.

Recall that if a real entire function $\varphi$ belongs to the
Laguerre-P\'olya class $\mathcal{LP}$ then satisfies the Laguerre
type inequalities
\begin{equation}\label{lag1}\left[\varphi^{(n)}(x)\right]^2-\varphi^{(n-1)}(x)\varphi^{(n+1)}(x)\geq0, \ \ n\in\{1,2,\dots\}\end{equation}
and
\begin{equation}\label{lag2}
\sum_{j=0}^{2k}\frac{(-1)^{j+k}}{(2k)!}\left(\begin{array}{c}2k\\
j\end{array}\right)\varphi^{(n+j)}(x)\varphi^{(n+2k-j)}(x)\geq0,\ \
\ n,k\in\{0,1,\dots\}.
\end{equation}
The inequality \eqref{lag1} is due to Skovgaard \cite{skov}, while
\eqref{lag2} for $n=0$ has been proved first by Jensen \cite{jensen}
and for $n\in\{0,1,\dots\}$ by Patrick \cite{patrick}.

Now, we list some other necessary conditions for a real entire
function to belong to $\mathcal{LP}$. For the inequality \eqref{t1}
see \cite{craven1,csordas}, the inequality \eqref{t3} is due to
Dimitrov \cite[Theorem 1]{dimitrov}, while \eqref{t4} can be found
in Patrick's paper \cite[Theorem 2]{patrick}. We note that when
$k=1$ the inequality \eqref{t4} reduces to \eqref{t1}.

\begin{theorem}\label{thlagpo}
Let $\{b_n(t)\}_{n\geq0}$ be a sequence of real functions which for
certain values of $t$ have a generating function of the type
\eqref{entire} and suppose that the real entire function $\varphi,$
defined by \eqref{entire}, is in the class $\mathcal{LP}$. Then for
those values of $t$ the following Tur\'an type and higher order
Tur\'an type inequalities hold
\begin{equation}\label{t1}
b_n^2(t)-b_{n-1}(t)b_{n+1}(t)\geq0,
\end{equation}
\begin{align}\label{t3}
4&\left[b_n^2(t)-b_{n-1}(t)b_{n+1}(t)\right]\left[b_{n+1}^2(t)-b_{n}(t)b_{n+2}(t)\right]\nonumber\\&-\left[b_n(t)b_{n+1}(t)-b_{n-1}(t)b_{n+2}(t)\right]^2\geq0,
\end{align}
\begin{equation}\label{t4}
\sum_{j=0}^{2k}\frac{(-1)^{j+k}}{(2k)!}\left(\begin{array}{c}2k\\
j\end{array}\right)b_{n+j}(t)b_{n+2k-j}(t)\geq0,\ \ \
n,k\in\{0,1,\dots\},
\end{equation}
where in \eqref{t1} and \eqref{t3} it is assumed that
$n\in\{1,2,\dots\}.$
\end{theorem}

\subsection{Higher order Tur\'an type inequalities for Bessel
functions of natural order} Now, we are going to point out other
results on Bessel functions of the first kind. Some of these results
are known, however we have included here for the sake of
completeness. We start with the formula
\begin{equation}\label{b1}
J_0\left(\sqrt{t^2-2xt}\right)=\sum_{n\geq0}J_n(t)\frac{x^n}{n!},\ \
\ t\in\mathbb{R},
\end{equation}
which is actually a special case of Bessel's formula (proved later
also by von Lommel) \cite[p. 140]{watson}
$$(z+h)^{-\nu/2}J_{\nu}\left(\sqrt{z+h}\right)=\sum_{n\geq0}\frac{(-h/2)^n}{n!}z^{-(\nu+n)/2}J_{\nu+n}\left(\sqrt{z}\right)$$
rewritten in the form
\begin{equation}\label{b2}
t^{\nu}(t^2-2xt)^{-\nu/2}J_{\nu}\left(\sqrt{t^2-2xt}\right)=\sum_{n\geq0}J_{\nu+n}(t)\frac{x^n}{n!},\
\ \ t\in\mathbb{R}.
\end{equation}
Since the Bessel function $J_0$ belongs to the Laguerre-P\'olya
class $\mathcal{LP}$ (see \cite[p. 123]{polya}), by using the
particular Bessel formula \eqref{b1}, we can see easily that the
Bessel functions of the first kind $J_n$ satisfies the inequalities
\eqref{t1}, \eqref{t3} and \eqref{t4}, that is, for all
$t\in\mathbb{R}$ we have the inequalities
   \begin{equation} \label{bt1}
      \Delta_n(t)\geq0,\ \ \ n\in\{1,2,\dots\},
   \end{equation}
   \begin{equation}\label{bt2}
      4\left[\Delta_n(t)\right]\cdot\left[\Delta_{n+1}(t)\right]-\left[J_n(t)J_{n+1}(t)-J_{n-1}(t)J_{n+2}(t)\right]^2\geq0,\
             \ \ n\in\{1,2,\dots\},
   \end{equation}
   \begin{equation} \label{bt3}
      \sum_{j=0}^{2k}\frac{(-1)^{j+k}}{(2k)!}\left(\begin{array}{c}2k\\
      j\end{array}\right)J_{n+j}(t)J_{n+2k-j}(t)\geq0,\ \ \  n,k\in\{0,1,\dots\}.
   \end{equation}
The Tur\'an type inequality \eqref{bt1} is well-known, it has been
discussed in details in the previous section. The inequality
\eqref{bt3} has been deduced by Patrick \cite{patrick}, and
\eqref{bt2} seems to be new. We note that, as it was pointed out by
Skovgaard \cite{skov} in 1954 and later by Patrick \cite{patrick} in
1973, the inequalities \eqref{bt1} and \eqref{bt3} are also
satisfied by the first derivative $J_n'$ of the Bessel function of
the first kind $J_n.$ More precisely, since the series
$$\sum_{n\geq0}J_n'(t)\frac{x^n}{n!}$$
is uniformly convergent for $t\in\mathbb{R},$ differentiating both
sides of \eqref{b1} with respect to $t$ and by using the relation
$J_0'(z)=-J_1(z)$ we get
$$(x-t)(t^2-2xt)^{-1/2}J_1\left(\sqrt{t^2-2xt}\right)=\sum_{n\geq0}J_n'(t)\frac{x^n}{n!}.$$
Now taking into account that $J_0\in\mathcal{LP}$ and the
Laguerre-P\'olya class is closed under differentiation, the
left-hand side of the above relation also belongs to the
Laguerre-P\'olya class $\mathcal{LP}$ and consequently (in view of
Theorem \ref{thlagpo}) for all $t\in\mathbb{R}$ we have the
following inequalities
\begin{equation} \label{btd1}
      [J_n'(t)]^2-J_{n-1}'(t)J_{n+1}'(t)\geq0,
   \end{equation}
   \begin{align}\label{btd2}
      4\left[[J_n'(t)]^2-J_{n-1}'(t)J_{n+1}'(t)\right]&\left[[J_{n+1}'(t)]^2-J_{n}'(t)J_{n+2}'(t)\right]\nonumber\\&-
      \left[J_n'(t)J_{n+1}'(t)-J_{n-1}'(t)J_{n+2}'(t)\right]^2\geq0,
   \end{align}
   \begin{equation} \label{btd3}
      \sum_{j=0}^{2k}\frac{(-1)^{j+k}}{(2k)!}\left(\begin{array}{c}2k\\
      j\end{array}\right)J_{n+j}'(t)J_{n+2k-j}'(t)\geq0,\ \ \  n,k\in\{0,1,\dots\}.
   \end{equation}
The inequality \eqref{btd1} is due to Skovgaard \cite{skov} and
holds for all $n\in\{1,2,\dots\},$ while the inequality \eqref{btd3}
is due to Patrick \cite{patrick}. The inequality \eqref{btd2}
appears to be new and holds for all $n\in\{1,2,\dots\}.$ Moreover,
as we can see below, it can be shown that in fact the above
inequalities for $J_n$ and $J_n'$ can be extended to real
parameters.

\subsection{Higher order Tur\'an type inequalities for Bessel
functions of real order} In order to extend the inequalities
\eqref{bt1}, \eqref{bt2}, \eqref{bt3}, \eqref{btd1}, \eqref{btd2}
and \eqref{btd3} to the Bessel function of the first kind $J_{\nu}$
and its derivative $J_{\nu}',$ where $\nu>-1,$ we shall use the
following result of Krasikov \cite{krasikov} and of Dimitrov and Ben
Cheikh \cite{dimit}.

\begin{lemma}
The function $\mathcal{J}_{\nu}:\mathbb{R}\to(-\infty,1],$ defined
by $\mathcal{J}_{\nu}(x)=2^{\nu}\Gamma(\nu+1)x^{-\nu}J_{\nu}(x),$
belongs to the Laguerre-P\'olya class $\mathcal{LP}$ when $\nu>-1.$
\end{lemma}

Note that Krasikov's \cite{krasikov} argument is simply and use only
the infinite product representation \eqref{weier}, by noticing that
the exponential factors in \eqref{weier} are canceled because of the
symmetry of the zeros $j_{\nu,n}$ with respect to the origin. The
proof of the above lemma given by Dimitrov and Ben Cheikh
\cite{dimit} use the following well-known result of Jensen: an
entire function $\varphi,$ defined by \eqref{entire} belongs to the
Laguerre-P\'olya class $\mathcal{LP}$ if and only if the associated
Jensen polynomials $g_n(\varphi;x),$ defined by
$$g_n(\varphi;x)=\sum_{k=0}^n\left(\begin{array}{c}n\\
k\end{array}\right)b_kx^k,$$ have only real zeros for all
$n\in\{0,1,\dots\}.$ It is also important to note that in fact
Dimitrov and Ben Cheikh \cite{dimit} by using the above lemma
reobtained von Lommel's celebrated result: the zeros of the Bessel
function $J_{\nu}$ are all real when $\nu>-1.$ Moreover, they proved
that the Jensen polynomials associated with the Bessel function
$J_{\nu}$, properly normalized, are exactly the Laguerre
polynomials. More precisely, in \cite{dimit} it was shown that the
only Jensen polynomials associated with an entire function in the
Laguerre-P\'olya class $\mathcal{LP}$ that are orthogonal are the
Laguerre polynomials.

Now, we are in the position to state and prove the following inequalities.

\begin{theorem}
Let $\nu>-1$ and $t\in\mathbb{R}.$ Then the following Laguerre type
inequalities
\begin{equation}\label{llag1}
\left[\mathcal{J}_{\nu}^{(n)}(t)\right]^2-\mathcal{J}_{\nu}^{(n-1)}(t)\mathcal{J}_{\nu}^{(n+1)}(t)\geq0,
\ \ n\in\{1,2,\dots\},
\end{equation}
\begin{equation}\label{llag2}
\sum_{j=0}^{2k}\frac{(-1)^{j+k}}{(2k)!}\left(\begin{array}{c}2k\\
j\end{array}\right)\mathcal{J}_{\nu}^{(n+j)}(t)\mathcal{J}_{\nu}^{(n+2k-j)}(t)\geq0,\
\ \ n,k\in\{0,1,\dots\}
\end{equation}
 and Tur\'an type inequalities
\begin{equation}\label{btn1}
J_{\nu+1}^2(t)-J_{\nu}(t)J_{\nu+2}(t)\geq0,
\end{equation}
\begin{align}\label{btn2}
4\left[J_{\nu+1}^2(t)-J_{\nu}(t)J_{\nu+2}(t)\right]&\left[J_{\nu+2}^2(t)-J_{\nu+1}(t)J_{\nu+3}(t)\right]\nonumber\\
&-\left[J_{\nu+1}(t)J_{\nu+2}(t)-J_{\nu}(t)J_{\nu+3}(t)\right]^2\geq0,
\end{align}
\begin{equation}\label{btn3}
\sum_{j=0}^{2k}\frac{(-1)^{j+k}}{(2k)!}\left(\begin{array}{c}2k\\
      j\end{array}\right)J_{\nu+j}(t)J_{\nu+2k-j}(t)\geq0,\ \ \  k\in\{0,1,\dots\},
\end{equation}
\begin{equation}\label{btnd1}
\left[J_{\nu+1}'(t)\right]^2-J_{\nu}'(t)J_{\nu+2}'(t)\geq0,
\end{equation}
\begin{align}\label{btnd2}
4\left[\left[J_{\nu+1}'(t)\right]^2-J_{\nu}'(t)J_{\nu+2}'(t)\right]&\left[\left[J_{\nu+2}'(t)\right]^2-J_{\nu+1}'(t)J_{\nu+3}'(t)\right]\nonumber\\&
-\left[J_{\nu+1}'(t)J_{\nu+2}'(t)-J_{\nu}'(t)J_{\nu+3}'(t)\right]^2\geq0,
\end{align}
\begin{equation}\label{btnd3}
\sum_{j=0}^{2k}\frac{(-1)^{j+k}}{(2k)!}\left(\begin{array}{c}2k\\
      j\end{array}\right)J_{\nu+j}'(t)J_{\nu+2k-j}'(t)\geq0,\ \ \  k\in\{0,1,\dots\},
\end{equation}
are valid.
\end{theorem}

\begin{proof}
Inequalities \eqref{llag1} and \eqref{llag2} follow immediately from
Lemma and the Laguerre type inequalities \eqref{lag1} and
\eqref{lag2}. By using Lemma, the left-hand side of
\eqref{b2} belongs to the Laguerre-P\'olya class $\mathcal{LP}$ and
consequently Theorem \ref{thlagpo} guarantees that the inequalities
\eqref{bt1}, \eqref{bt2} and \eqref{bt3} hold true for $\nu+n$
instead of $n.$ Now, choosing in the modified form of \eqref{bt1},
\eqref{bt2} the value $n=1$ and in the modified form of \eqref{bt3}
the value $n=0,$ we get the inequalities \eqref{btn1}, \eqref{btn2}
and \eqref{btn3}.

By using the idea of Skovgaard, since the series
$$\sum_{n\geq0}J_{\nu+n}'(t)\frac{x^n}{n!}$$
is uniformly convergent for $t\in\mathbb{R},$ differentiating both
sides of \eqref{b2} with respect to $t$ we get
$$\left[t^{\nu}(t^2-2xt)^{-\nu/2}J_{\nu}\left(\sqrt{t^2-2xt}\right)\right]'=\sum_{n\geq0}J_{\nu+n}'(t)\frac{x^n}{n!},\
\ \ t\in\mathbb{R}.$$ Now, because of Lemma the left-hand
side of \eqref{b2} belongs to $\mathcal{LP},$ the left-hand side of
the above relation belongs also to $\mathcal{LP}.$ This together
with Theorem \ref{thlagpo} imply that the Tur\'an type inequalities
\eqref{btd1}, \eqref{btd2} and \eqref{btd3} hold true for $\nu+n$
instead of $n.$ Similarly, as above, by choosing in the modified
form of \eqref{btd1}, \eqref{btd2} the value $n=1$ and in the
modified form of \eqref{btd3} the value $n=0,$ we get the Tur\'an
type inequalities \eqref{btnd1}, \eqref{btnd2} and \eqref{btnd3}.
\end{proof}

\begin{remark}
{\em We note that when $k=1$ the inequality \eqref{btn3} reduces to
\eqref{btn1}, while the inequality \eqref{btnd3} becomes
\eqref{btnd1}. The inequality \eqref{btn1} is exactly \eqref{turan1}
for $\nu>0$ and the inequality \eqref{btn3} bears resemblance of
\eqref{salamt}. We also note here that the inequality \eqref{llag2}
for $n=0$ was one of the crucial facts in Krasikov's \cite{krasikov}
study in order to obtain sharp uniform bounds for the Bessel
functions of the first kind $J_{\nu}.$}
\end{remark}

It is worth to mention here that the Tur\'an type inequality
\eqref{turan2} can be rewritten in terms of $\mathcal{J}_{\nu}$ as
$$\mathcal{J}_{\nu}^2(x)-\mathcal{J}_{\nu-1}(x)\mathcal{J}_{\nu+1}(x)\geq0,$$
where $\nu>0$ and $x\in\mathbb{R}.$ A counterpart of the above
Tur\'an type inequality was derived by Sz\'asz \cite{szasz2} in 1951
(and also by Joshi and Bissu \cite{JB} in 1991) as follows
$$(\nu+1)\mathcal{J}_{\nu}^2(x)-\nu\mathcal{J}_{\nu-1}(x)\mathcal{J}_{\nu+1}(x)\leq1,$$
which also holds for all $\nu>0$ and $x\in\mathbb{R}.$ These
inequalities together show that the quantity
$$(\nu+1)\mathcal{J}_{\nu}^2(x)-\nu\mathcal{J}_{\nu-1}(x)\mathcal{J}_{\nu+1}(x)$$
belongs to the interval $[0,1]$ for all $x\in\mathbb{R}$ and
$\nu>0.$ Moreover, by means of the differentiation formulae
$$\mathcal{J}_{\nu}'(x)=-\frac{x}{2(\nu+1)}\mathcal{J}_{\nu+1}(x),$$
$$\mathcal{J}_{\nu}''(x)=-\frac{1}{2(\nu+1)}\mathcal{J}_{\nu+1}(x)+\frac{x^2}{4(\nu+1)(\nu+2)}\mathcal{J}_{\nu+2}(x),$$
the Laguerre type inequality \eqref{llag1} for $n=1$ becomes
$$\frac{x^2}{2(\nu+1)}\mathcal{J}_{\nu+1}^2(x)+\mathcal{J}_{\nu}(x)\mathcal{J}_{\nu+1}(x)\geq
\frac{x^2}{2(\nu+2)}\mathcal{J}_{\nu}(x)\mathcal{J}_{\nu+2}(x)$$
which is equivalent to
$$J_{\nu+1}^2(x)-J_{\nu}(x)J_{\nu+2}(x)\geq -\frac{1}{x}J_{\nu}(x)J_{\nu+1}(x), \ \ \ x\neq0,\ \nu>-1.$$

\section{\bf Integral representations for second kind Neumann type series of Bessel functions}
\setcounter{equation}{0}

The applications of the Neumann series of Bessel functions
   \[ \mathfrak N_\nu(x) = \sum_{n\ge 1} \alpha_n J_{\nu+n}(x) \]
in science, engineering and technology has been scattered widely in
the literature, see for example \cite{PS} and the references
therein. Recently, Pog\'any and S\"uli \cite{PS} deduced closed
integral expression for $\mathfrak N_\nu(x)$ and with this initiated
considerable interest to closed integral expressions for Neumann
series of another kind of Bessel functions, see e.g. the results by
Baricz et al. \cite{BJP1, BJP2}.

In this section our aim is to deduce a closed integral
representation formula for the second kind Neumann type series of
Bessel functions defined in the form
   \begin{equation} \label{N1}
      \mathfrak G_{\mu, \nu}^{a,b}(x) := \sum_{n \ge 1} \theta_n J_{\mu+an}(x)J_{\nu+bn}(x), \qquad \mu, \nu, a, b\in \mathbb R\, .
   \end{equation}
This is motivated by the fact that $\mathfrak G_{\nu, \nu}^{2,2}(x)$
stands for the right-hand side series in von Lommel's expression
\eqref{lommel} and for the Al--Salam series \cite{salam} appearing
in \eqref{salam}, while $\mathfrak G_{\nu, \nu}^{1,1}(x)$ covers the
series considered by Thiruvenkatachar and Nanjundiah \cite{thiru} in
\eqref{turan4}; finally, Neumann type finite sum close to $\mathfrak
G_{\,\,\nu, \nu}^{1,-1}(x)$ appears in \eqref{salam}, and another
finite sum of type $\mathfrak G_{n, n+2k}^{\,\,1,-1}(x)$ takes place
in \eqref{bt3}.

In order to obtain the integral representation formula for
\eqref{N1} we shall use the main idea from \cite{PS}, that is, the
Laplace integral representation of the associated Dirichlet series
will be the main tool here. Thus, we take $x\in \mathbb R_+$ and
assume in the sequel that the behavior of $\{\theta_n\}_{n\ge 1}$
ensures the convergence of the series \eqref{N1} over $\mathbb R_+$.
Throughout $[a]$ and $\{a\}= a-[a]$ denote the integer and
fractional part some real $a$ respectively and $\chi_S$ will stands
for the characteristic function of the set $S\subset \mathbb{R}$.
Moreover, let us consider the real-valued function $x\mapsto a_x =
a(x)$ and suppose that $a\in {\rm C}^1[k, m]$, $k,m \in \mathbb Z$,
$k<m$. The classical Euler--Maclaurin summation formula states that
   \[ \sum_{j=k}^m a_j = \int_k^m a(x){\rm d}x + \frac12\left(a_k+a_m\right) + \int_k^m \left( x-[x]-\dfrac 12\right)a'(x) {\rm d}x.\]
The following condensed form of the Euler--Maclaurin formula holds \cite{PS}:
   \begin{equation} \label{N2}
      \sum_{j=k+1}^m a_j = \int_k^m \big( a(x)+\{ x\}a'(x) \big) {\rm d}x = \int_k^m \mathfrak d_x a(x)\, {\rm d}x\, ,
   \end{equation}
where
   \[  \mathfrak d_x := 1 + \{ x \}\frac{\rm d}{{\rm d} x}.\]
Also, we need in the sequel a tool to estimate Bessel functions of the first kind. Landau \cite{LL} gave in a sense
the best possible uniform bound for the first kind Bessel function $J_\nu$ with respect to $x, \nu>0$:
   \begin{equation} \label{X6}
      |J_\nu(x)| \le b_L\, \nu^{-1/3}, \qquad b_L = \sqrt[3]{2}\sup_{x\in \mathbb R_+} {\rm Ai}(x),
   \end{equation}
where $\text{Ai}(\cdot)$ stands for the familiar Airy function
\cite[p. 447]{abra}. In fact Krasikov \cite{krasikov} pointed out
that this inequality is sharp only in the transition region, i.e.
for $x$ around $j_{\nu,1}$, the first positive zero of $J_\nu.$ For
further reading about and detailed discussions consult \cite{O} and
\cite{PS1} where another fashion upper bounds are given.

\begin{theorem}\label{thp}
Let $\theta\in {\rm C}^1(\mathbb R_+),$ $\theta\big|_{\mathbb N} =
\{\theta_n\}_{n\ge 1}$ and assume that series $\sum_{n \ge
1}\theta_n\, n^{-2/3}$ is absolutely convergent. Then, for all
$x,a,b>0$
 such that
   \begin{equation} \label{N3}
      x \in \left( 0, 2\min\left\{1, \dfrac1{e} \left( \dfrac{a^ab^b}{\rho_{\,\,\mathfrak G}^{a,b}}\right)^{1/(a+b)}\right\}\right) =
                      \mathcal I_{\mathfrak G}, \qquad \min(\mu+a, \nu+b)>0 ,
   \end{equation}
where
   \begin{equation} \label{N31}
      \rho_{\,\,\mathfrak G}^{a,b} = \limsup_{n\to \infty}\frac{|\theta_n|^{1/n}}{n^{a+b}}\, ,
   \end{equation}
we have that
   \begin{align} \label{N5}
      \mathfrak G_{\mu, \nu}^{a,b}(x) &= -\int_1^\infty \int_0^{[u]} \frac{\partial}{\partial u}
                            \left(\,\Gamma(\mu+au+1)\Gamma(\nu+bu+1)\, J_{\mu+au}(x)J_{\nu+bu}(x)\right) \nonumber \\
                 & \qquad \times \mathfrak d_v\left( \frac{\theta(v)}{\Gamma(\mu+av)\Gamma(\nu+bv)}\right){\rm d}u \,{\rm d}v\, .
   \end{align}
\end{theorem}

\begin{proof} Applying Landau's bound \eqref{X6} we estimate Bessel functions of the first kind such that consists $\mathfrak G_{\mu, \nu}^{a,b}(x)$
   \[ \big| \mathfrak G_{\mu, \nu}^{a,b}(x)\big| \le b_L^2 \sum_{n \ge 1}\dfrac{|\theta_n|}{\sqrt[3]{(\mu+an)(\nu+bn)} }
            \sim \dfrac{b_L^2}{\sqrt[3]{ab}} \sum_{n \ge 1} \dfrac{|\theta_n|}{n^{2/3}}\, ,\]
hence $\mathfrak G_{\mu, \nu}^{a,b}(x)$ absolutely and uniformly
converges for $x>0$. Now, by substituting the integral
representation formula \cite[p. 48]{watson}
   \[ J_\nu(x) = \frac{2\left(x/2\right)^\nu}{\sqrt{\pi}\, \Gamma(\nu+1/2)}\int_0^1 \cos(xt)(1-t^2)^{\nu-\frac{1}{2}}{\rm d}t,\qquad
      x\in \mathbb R,\,\nu> -1/2, \]
into \eqref{N1} we conclude
   \begin{equation} \label{N6}
      \mathfrak G_{\mu, \nu}^{a,b}(x) = \dfrac4\pi \left( \frac x2\right)^{\mu+\nu} \int_0^1\int_0^1 \cos(xt)\cos(xs)
                   (1-t^2)^{\mu-\frac{1}{2}}(1-s^2)^{\nu-\frac{1}{2}}\, \mathcal D_\theta(t,s)\, {\rm d}t{\rm d}s\, ,
   \end{equation}
where the Dirichlet series'
   \[ \mathcal D_\theta(t,s) = \sum_{n \ge 1}\dfrac{\theta_n \left(\left(x/2\right)^{a+b}(1-t^2)^a(1-s^2)^b\right)^n}{\Gamma(\mu+an)\Gamma(\nu+bn)} \]
$x$--domain of convergence is our first goal. To express Dirichlet
series in Laplace integral form it is necessary to have positive
Dirichlet parameter, that is,
   \[ -\ln \big(x/2\big)^{a+b}(1-t^2)^a(1-s^2)^b >0,\]
which holds for all $|x|<2$ when $a+b>0.$ Also $\mathcal
D_\theta(t,s)$ is equiconvergent to the auxiliary power series
   \[ \sum_{n \ge 1} \dfrac{\theta_n}{n^{(a+b)n}}\,\left(\dfrac{({e}x)^{a+b}(1-t^2)^a(1-s^2)^b}{2^{a+b}a^ab^b} \right)^n\]
with radius of convergence
   \[ \rho_{\,\,\mathfrak G}^{a,b} =  \left(\limsup_{n \to \infty}\dfrac{|\theta_n|^{1/n}}{n^{a+b}}\right)^{-1}\, .\]
This yields the convergence region $\mathcal I_{\mathfrak G}$
described in \eqref{N3}.

Next, the Laplace integral expression (cf. \cite{K}) for Dirichlet
series $\mathcal D_\theta(t,s)$ becomes
   \begin{align} \label{N7}
      \mathcal D_\theta(t,s) &= \sum_{n \ge 1}\dfrac{\theta_n }{\Gamma(\mu+an)\Gamma(\nu+bn)}\,
               \exp \left\{ - n \ln \dfrac{2^{a+b}}{x^{a+b}(1-t^2)^{a}(1-s^2)^{b}}\right\} \nonumber \\
            &= \ln \dfrac{2^{a+b}}{x^{a+b}(1-t^2)^{a}(1-s^2)^{b}} \int_0^\infty
               \left(\dfrac{x^{a+b}(1-t^2)^{a}(1-s^2)^{b}}{2^{a+b}}\right)^u \nonumber \\
            & \qquad \times  \left(\sum_{n=1}^{[u]} \frac{\theta_n}{\Gamma(\nu+an)\Gamma(\mu+bn)}\right) {\rm d}u \nonumber \\
            &= \ln \dfrac{2^{a+b}}{x^{a+b}(1-t^2)^{a}(1-s^2)^{b}} \int_0^\infty
               \left(\dfrac{x^{a+b}(1-t^2)^{a}(1-s^2)^{b}}{2^{a+b}}\right)^u  \nonumber \\
            &\qquad \times \int_0^{[u]} \mathfrak d_v \left( \dfrac{\theta(v)\, {\rm d}v}{\Gamma(\nu+an)\Gamma(\mu+bn)}\right)\, {\rm d}u\, ,
   \end{align}
where the last equality we deduced by virtue of condensed
Euler--Maclaurin summation formula \eqref{N2}. The last formula in
conjunction with \eqref{N6} gives
   \[ \mathfrak G_{\mu, \nu}^{a,b}(x) = \dfrac4\pi \left( \frac x2\right)^{\mu+\nu} \int_0^\infty\int_0^{[u]}
                \mathcal J_{t,s}(u)\, \mathfrak d_v \left( \dfrac{\theta(v)}{\Gamma(\nu+an)\Gamma(\mu+bn)}\right)\, {\rm d}u{\rm d}v\, ,\]
where
   \begin{align*}
      \mathcal J_{t,s}(u) &= - \left(\dfrac x2\right)^{(a+b)u}\int_0^1\int_0^1 \cos(xt)\cos(xs)\,
                (1-t^2)^{\mu+au-\frac{1}{2}}(1-s^2)^{\nu+bv-\frac{1}{2}} \nonumber \\
            &\qquad \times \ln \left(\dfrac{x^{a+b}(1-t^2)^{a}(1-s^2)^{b}}{2^{a+b}}\right) {\rm d}t{\rm d}s\, \, .
   \end{align*}
Because
   \[ \int \mathcal J_{t,s}(u)\, {\rm d}u = - \left(\dfrac x2\right)^{(a+b)u} \mathcal I_{\mu,a}(u)\mathcal I_{\nu,b}(u)\, ,\]
where
   \begin{align*}\mathcal I_{\mu,a}(u)& = \int_0^1 \cos(xt) (1-t)^{\mu+au-\frac{1}{2}}{\rm
   d}t\\& = \dfrac{\sqrt{\pi}}2 \left( \dfrac x2\right)^{-\mu-au}
   \Gamma(\mu+au+1)J_{\mu+au}(x)\end{align*}
there holds
   $$
      \mathcal J_{t,s}(u) = - \dfrac\pi4 \left(\dfrac x2\right)^{-\mu-\nu}\, \dfrac{\partial}{\partial u}
               \left(\Gamma(\mu+au+1)\Gamma(\nu+bu+1)J_{\mu+au}(x)J_{\nu+bu}(x)\right)\, .
   $$
Now, straightforward transformations yield the asserted formula \eqref{N5}.
\end{proof}

In the preliminary part of this section we mentioned the equalities
by von Lommel, Thiruvenkatachar and Nanjundiah and Al--Salam, where
previous integral expression's special cases play important roles.
By particular choices of $a$ and $b$ we have the following cases.

\begin{corollary} If $\nu>2$ and
   \[ x \in \left( 0, 2\min\left\{1, 2\left({e}^4{\rho_{\,\,\mathfrak G}^{2,2}}\right)^{-\frac{1}{4}}\right\}\right),\]
then we have that
   \begin{align} \label{N9}
      \mathfrak G_{\nu-2, \nu-2}^{\quad2,2}(x) &= - \int_1^\infty\int_0^{[u]} \dfrac{\partial}{\partial u}\left(
                \Gamma^2(\nu-1+2u) J_{\nu-2+2u}^2(x)\right)\nonumber \\
             & \qquad \times \mathfrak d_v \left( \dfrac{(\nu-2+2v)\Gamma(-1+v)\Gamma(\nu-1+v)}
                {\Gamma^2(\nu-2-2v)\Gamma(v)\Gamma(\nu+v)}\right)\, {\rm d}u{\rm d}v\, .
   \end{align}
Moreover, for $\nu>0$ and
   \[ x \in \left( 0, 2\min\left\{1, \left({e}^2 \rho_{\,\,\mathfrak G}^{1,1}\right)^{-\frac{1}{2}}\,\right\}\right)\]
there holds true
   \begin{align} \label{N10}
      \mathfrak G_{\nu+1, \nu+1}^{\quad1,1}(x) &= - \int_1^\infty\int_0^{[u]} \dfrac{\partial}{\partial u}\left(\Gamma^2(\nu+u+2)
                J_{\nu+u+2}^2(x)\right)\nonumber \\
                & \qquad \times \mathfrak d_v \left(\dfrac{(\nu+v)^{-1}(\nu+v+2)^{-1}}{\Gamma^2(\nu+v+1)}\right)\, {\rm d}u{\rm d}v.
   \end{align}
\end{corollary}

We note that the second kind Neumann series in von Lommel's formula
\eqref{lommel} possesses divergent auxiliary series $\sum_{n \ge
1}\theta_n \cdot n^{-2/3}$, therefore it is not covered by Theorem
\ref{thp}. Also, Al-Salam's Neumann series, taking place in
\eqref{salam}, converges only when $m<1/3$, so the integral
representation \eqref{N9}. The auxiliary series associated with the
Neumann series by Thiruvenkatachar and Nanjundiah in \eqref{turan4}
is equiconvergent to the Riemannian $\zeta(8/3)$, thus this case
meet Theorem \ref{thp}, see \eqref{N10}. Finally, the second kind
finite sums \eqref{salam} and \eqref{bt3} of Neumann type have to be
considered separately, since their summation require different
approach, being the parameter space of $\mathfrak G_{\mu,
\nu}^{a,b}(x)$ restricted by the required positivity of upper
parameters $a,b$ in Theorem \ref{thp}.

\section{\bf Tur\'an type inequalities for Bessel functions of the
second kind and related problems} \setcounter{equation}{0}

\subsection{Tur\'an type inequalities for Bessel functions of the
second kind} The Bessel function of the second kind $Y_{\nu}$ is
defined by \cite[p. 64]{watson}
$$Y_{\nu}(x)=\frac{J_{\nu}(x)\cos(\nu\pi)-J_{-\nu}(x)}{\sin(\nu\pi)},$$
where the right-hand side of this equation is replaced by its
limiting value if $\nu$ is an integer or zero. Recently, Baricz
\cite{bams} posed the following open problem concerning the Bessel
function of the second kind: Is it true that the Tur\'an type
inequality
\begin{equation}\label{turany0}
Y_{\nu}^2(x)-Y_{\nu-1}(x)Y_{\nu+1}(x)>\frac{1}{1-\nu}Y_{\nu}^2(x)
\end{equation}
holds true for all $\nu\in\mathbb{R}$ and $x>0$? For those values of
$\nu$ for which \eqref{turany0} holds the constant $1/(1-\nu)$ is
the best possible?

In this section our aim is to give the solution to the above open
problem. Observe that the Bessel function of the second kind
$Y_{\nu}$ has the same recurrence relations as the Bessel function
of the first kind $J_{\nu}$ and hence by using the idea from
\cite{thiru} we may obtain the following result.

\begin{theorem}\label{thy}
If $\nu>1$ and $x>x_{\nu},$ where $x_{\nu}\leq\nu$ is the unique
positive root of the equation
$$Y_{\nu}^2(x)-Y_{\nu-1}(x)Y_{\nu+1}(x)=0,$$ then
the Tur\'an type inequality
\begin{equation}\label{turany1}
Y_{\nu}^2(x)-Y_{\nu-1}(x)Y_{\nu+1}(x)>0
\end{equation}
is valid. Moreover, \eqref{turany1} is reversed for $0<x<x_{\nu}.$
\end{theorem}

\begin{proof}
Let us denote
$$\overline{\Delta}_{\nu}(x)=Y_{\nu}^2(x)-Y_{\nu-1}(x)Y_{\nu+1}(x)$$
and consider the function $\Psi_{\nu}:(0,\infty)\to\mathbb{R},$
defined by $\Psi_{\nu}(x)=x^2\overline{\Delta}_{\nu}(x).$ In view of
the recurrence relations \cite[p. 66]{watson}
\begin{equation}\label{recy1}xY_{\nu}'(x)+\nu Y_{\nu}(x)=xY_{\nu-1}(x)\end{equation} and
\begin{equation}\label{recy2}xY_{\nu}'(x)-\nu Y_{\nu}(x)=-xY_{\nu+1}(x),\end{equation}
we get
\begin{align}\label{yroot}\Psi_{\nu}(x)&=x^2Y_{\nu}^2(x)-\left[\nu
Y_{\nu}(x)+xY_{\nu}'(x)\right]\left[\nu
Y_{\nu}(x)-xY_{\nu}'(x)\right]\nonumber\\&
=(x^2-\nu^2)Y_{\nu}^2(x)+x^2\left[Y_{\nu}'(x)\right]^2\end{align}
and then
$$\Psi_{\nu}'(x)=2x\left[Y_{\nu}(x)\right]^2+2Y_{\nu}'(x)\left[x^2Y_{\nu}''(x)+xY_{\nu}'(x)+(x^2-\nu^2)Y_{\nu}(x)\right].$$
Since the Bessel function of the second kind $Y_{\nu}$ satisfies
\eqref{eqbessel}, we obtain
$$x^2Y_{\nu}''(x)+xY_{\nu}'(x)+(x^2-\nu^2)Y_{\nu}(x)=0,$$
which in turn implies that
$\Psi_{\nu}'(x)=2x\left[Y_{\nu}(x)\right]^2,$ and then the function
$\Psi_{\nu}$ is increasing on $(0,\infty).$

Now recall the asymptotic formula \cite[p. 360]{abra}
$$Y_{\nu}(x)\sim -\frac{1}{\pi}\Gamma(\nu)\left(\frac{x}{2}\right)^{-\nu}$$
which holds when $\nu>0$ is fixed and $x\to 0.$ Then clearly when
$\nu>1$ is fixed and $x\to0$ we have
$$\Psi_{\nu}(x)\sim\frac{2^{2\nu}}{\pi^2}x^{2-2\nu}\left[\Gamma^2(\nu)-\Gamma(\nu-1)\Gamma(\nu+1)\right],$$
which in turn implies that $\Psi_{\nu}(x)$ tends to $-\infty$ as
$x\to 0.$ Here we used the fact that the Euler's gamma function
$\Gamma$ is logarithmically convex on $(0,\infty)$ and then for all
$\nu>1$ we have
$$\Gamma^2(\nu)<\Gamma(\nu-1)\Gamma(\nu+1).$$ Note that the above inequality can be realized also by means of the well-known formula
$\Gamma(\mu+1)=\mu\Gamma(\mu),$ $\mu>0.$ More precisely, since
$\Gamma(\mu)>0$ for all $\mu>0,$ we have
$$\Gamma^2(\nu)-\Gamma(\nu-1)\Gamma(\nu+1)=-\Gamma(\nu-1)\Gamma(\nu)<0.$$

On the other hand, it is known that when $\nu$ is fixed and
$x\to\infty$ we have \cite[p. 364]{abra}
$$Y_{\nu}(x)\sim\sqrt{\frac{2}{\pi x}}\sin\left(x-\nu\frac{\pi}{2}-\frac{\pi}{4}\right),$$
which yields
$$\Psi_{\nu}(x)\sim \frac{2}{\pi}x\left[\sin^2\left(x-\nu\frac{\pi}{2}-\frac{\pi}{4}\right)-
\sin\left(x-(\nu+1)\frac{\pi}{2}-\frac{\pi}{4}\right)\sin\left(x-(\nu-1)\frac{\pi}{2}-\frac{\pi}{4}\right)\right]=\frac{2}{\pi}x,$$
that is,  $\Psi_{\nu}(x)$ tends to $\infty$ as $x\to\infty.$

Summarizing, the function $\Psi_{\nu}$ is increasing, at zero tends
to $-\infty$ and at infinity tends to $\infty.$ Thus, there exists
an unique $x_{\nu}>0$ such that $\Psi_{\nu}(x)<0$ for $0<x<x_{\nu},$
$\Psi_{\nu}(x_{\nu})=0$ and $\Psi_{\nu}(x)>0$ for $x>x_{\nu}.$ With
this the proof of \eqref{turany1} is done. Finally, by using
\eqref{yroot} for all $x$ and $\nu$ real we have
$$\Psi_{\nu}(x)\geq (x^2-\nu^2)Y_{\nu}^2(x)$$ and then
$$0=\Psi_{\nu}(x_{\nu})\geq (x_{\nu}^2-\nu^2)Y_{\nu}^2(x_{\nu}),$$
which in turn implies that $x_{\nu}\leq\nu$ for all $\nu>0$ real.
\end{proof}

Observe that for $\nu>1$ and $x\geq x_{\nu}$ the inequality
\eqref{turany1} clearly implies \eqref{turany0}. Thus, although the
expression $1-Y_{\nu-1}(x)Y_{\nu+1}(x)/Y_{\nu}^2(x)$ tends to
$1/(1-\nu)$ if $\nu>1$ is fixed and $x$ tends to zero, because of
the next result, in \eqref{turany0} the constant $1/(1-\nu)$ is not
the best possible, even if was claimed in \cite{bams}.

We note also that by using the recurrence relation \cite[p.
66]{watson}
\begin{equation}\label{recy3}Y_{\nu-1}(x)+Y_{\nu+1}(x)=\frac{2\nu}{x}Y_{\nu}(x)\end{equation} we obtain
$$\nu Y_{\nu}(x)\left[Y_{\nu}(x)+Y_{\nu+2}(x)\right]=(\nu+1) Y_{\nu+1}(x)\left[Y_{\nu-1}(x)+Y_{\nu+1}(x)\right],$$
which in turn implies that
$$(\nu+1)\overline{\Delta}_{\nu}(x)-\nu\overline{\Delta}_{\nu+1}(x)=Y_{\nu}^2(x)+Y_{\nu+1}^2(x).$$
Consequently, the inequality
$$(\nu+1)\overline{\Delta}_{\nu}(x)-\nu\overline{\Delta}_{\nu+1}(x)>0$$
is valid for all admissible values of $\nu$ and $x.$ This inequality
in particular implies that for all $\nu>0$ we have
$\overline{\Delta}_{\nu}(x_{\nu+1})\geq0$ and
$\overline{\Delta}_{\nu+1}(x_{\nu})\leq0.$

By following the proof of Theorem \ref{thlak} we may prove the
following result, which completes Theorem \ref{thy}.
\begin{theorem}\label{thlaky}
For $\nu>0$ the relative maxima (denoted by $\overline{M}_{\nu,n}$)
of the function $x\mapsto \overline{\Delta}_{\nu}(x)$ occur at the
zeros of $Y_{\nu-1}$ and the relative minima (denoted by
$\overline{m}_{\nu,n}$) occur at the zeros of $Y_{\nu+1}.$ If
$\nu<0,$ the above statement is reversed. The values
$\overline{M}_{\nu,n}$ and $\overline{m}_{\nu,n}$ can be expressed
as
$$\overline{M}_{\nu,n}=\overline{\Delta}_{\nu}(y_{\nu-1,n})=Y_{\nu}^2(y_{\nu-1,n})$$
and
$$\overline{m}_{\nu,n}=\overline{\Delta}_{\nu}(y_{\nu+1,n})=Y_{\nu}^2(y_{\nu+1,n}),$$
where $y_{\nu,n}$ denotes the $n$th positive zero of the Bessel
function $Y_{\nu}.$ Consequently \eqref{turany1} holds true for all
$\nu\neq0$ and $x\geq y_{\nu-1,1}.$ If $\nu=0,$ then \eqref{turany1}
holds true for all $x\neq0.$
\end{theorem}

\begin{proof}
By using the recurrence relation \cite[p. 66]{watson}
$$Y_{\nu-1}(x)-Y_{\nu+1}(x)=2Y_{\nu}'(x),$$
and the formula \eqref{recy2} for $\nu-1$ instead of $\nu$ and
\eqref{recy1} for $\nu+1$ instead of $\nu,$ we obtain
$$\left[\overline{\Delta}_{\nu}(x)\right]'=\frac{2}{x}Y_{\nu-1}(x)Y_{\nu+1}(x)$$
and then
$$\left[\overline{\Delta}_{\nu}(x)\right]''=\frac{2}{x}Y_{\nu-1}'(x)Y_{\nu+1}(x)+
\frac{2}{x}Y_{\nu-1}(x)Y_{\nu+1}'(x)-\frac{2}{x^2}Y_{\nu-1}(x)Y_{\nu+1}(x).$$
Now, by using again \eqref{recy2} for $\nu-1$ instead of $\nu,$ and
\eqref{recy3}, we obtain
$$\left.\left[\overline{\Delta}_{\nu}(x)\right]'\right|_{x=y_{\nu-1,n}}=0\ \ \mbox{and} \ \
\left.\left[\overline{\Delta}_{\nu}(x)\right]''\right|_{x=y_{\nu-1,n}}=-\frac{4\nu}{y_{\nu-1,n}^2}Y_{\nu}^2(y_{\nu-1,n}).$$
Similarly, by using \eqref{recy1} for $\nu+1$ instead of $\nu,$ and
\eqref{recy3}, we get
$$\left.\left[\overline{\Delta}_{\nu}(x)\right]'\right|_{x=y_{\nu+1,n}}=0\ \ \mbox{and} \ \
\left.\left[\overline{\Delta}_{\nu}(x)\right]''\right|_{x=y_{\nu+1,n}}=\frac{4\nu}{y_{\nu+1,n}^2}Y_{\nu}^2(y_{\nu+1,n}).$$
These in turn imply that indeed for $\nu>0$ ($\nu<0$) the relative
maxima (minima) of the function $x\mapsto
\overline{\Delta}_{\nu}(x)$ occur at the zeros of $Y_{\nu-1}$ and
the relative minima (maxima) occur at the zeros of $Y_{\nu+1}.$ Now,
because the quantities $\overline{M}_{\nu,n}$ and
$\overline{m}_{\nu,n}$ are positive, \eqref{turany1} holds true for
all $\nu\neq0$ and $x\geq y_{\nu-1,1}.$ Finally, since
$Y_{-1}(x)=-Y_{1}(x),$ we have
$\overline{\Delta}_0(x)=Y_0^2(x)+Y_1^2(x),$ which is clearly
positive.
\end{proof}

\subsection{Lower bound for the logarithmic derivative of Bessel
functions of the second kind} In this subsection our aim is to
deduce new inequalities for the logarithmic derivative of $Y_{\nu}$
and for the quotient $Y_{\nu}/Y_{\nu-1}$ by using the Tur\'an type
inequality \eqref{turany1}. These results are the analogous of those
of subsection 2.7. Observe that by using Theorem \ref{thy} and
combining \eqref{turany1} with \eqref{recy1} and \eqref{recy2} we
immediately obtain the inequality
\begin{equation}\label{boundy1}
\left[xY_{\nu}'(x)\right]^2+(x^2-\nu^2)Y_{\nu}^2(x)\geq0,
\end{equation}
which holds for all $\nu>1$ and $x\geq x_{\nu}.$ Moreover, combining
\eqref{boundy1} with \eqref{recy1} we get
\begin{equation}\label{boundy2}
xY_{\nu}^2(x)-2\nu Y_{\nu-1}(x)Y_{\nu}(x)+xY_{\nu-1}^2(x)\geq0,
\end{equation}
where $\nu>1$ and $x\geq x_{\nu}.$ The inequality \eqref{boundy1}
implies the inequality
\begin{equation}\label{boundy3}
\left[\frac{xY_{\nu}'(x)}{Y_{\nu}(x)}\right]^2\geq \nu^2-x^2,
\end{equation}
where $\nu>1$ and $x\geq x_{\nu}$ such that $x\neq y_{\nu,n},$
$n\in\{1,2,\dots\},$ while the inequality \eqref{boundy2} implies
the following
\begin{equation}\label{boundy4} \frac{Y_{\nu}(x)}{Y_{\nu-1}(x)}\geq
\frac{\nu+\sqrt{\nu^2-x^2}}{x},
\end{equation}
where $\nu>1,$ $x_{\nu}\leq x\leq\nu$ such that $x\neq y_{\nu-1,n},$
$n\in\{1,2,\dots\}.$

Finally note that, according to Theorem \ref{thlaky} the conditions
$\nu>1$ and $x\geq x_{\nu}$ in inequalities \eqref{boundy1},
\eqref{boundy2} and \eqref{boundy3} can be replaced to $\nu\neq0$
and $x\geq y_{\nu-1,1},$ while in \eqref{boundy4} to $\nu>0$ and
$x\geq y_{\nu-1,1}.$

\subsection*{Acknowledgments} The research
of \'A. Baricz was supported in part by the J\'anos Bolyai Research
Scholarship of the Hungarian Academy of Sciences and in part by the
Romanian National Authority for Scientific Research CNCSIS-UEFISCSU,
project number PN-II-RU-PD\underline{ }388/2012. The work of this
author was initiated during his visit in April 2010 to the Indian
Institute of Technology Madras, Chennai, India. The visit was
supported by a partial travel grant from the Commission on
Development and Exchanges, International Mathematical Union. \'A.
Baricz acknowledges the hospitality from Professor Saminathan
Ponnusamy. Both of the authors are grateful to the referee for
his/her useful and constructive comments.

\end{document}